\documentclass[12pt,twoside]{article}

\usepackage{amssymb}
\usepackage{amsmath}
\usepackage{a4wide}
\usepackage{epic,eepic}
\usepackage{bbm}
\usepackage{amsthm}
\usepackage[all]{xy}

\newcommand{\id}{\operatorname{id}}
\newcommand{\Hom}{\operatorname{Hom}}
\newcommand{\lo}{\left(}
\newcommand{\po}{\right)}

\newcommand{\ulmod}{\operatorname{\underline{mod}}}

\newcommand{\End}{\operatorname{End}}

\newcommand{\md}{\operatorname{mod}}
\newcommand{\soc}{\operatorname{soc}}
\newcommand{\rad}{\operatorname{rad}}

\newcommand{\add}{\operatorname{add}}

\newcommand{\dualnosc}{\operatorname{D}}

\renewcommand{\tilde}{\widetilde}

\makeatletter

\newcommand{\circlearrowr}{
  \put(0,0){\oval(20,20)[r]}
  \bezier{50}(0,-10)(-4.14,-10)(-7.1,-7.1)
  \bezier{50}(0,10)(-4.14,10)(-7.1,7.1)
  \bezier{50}(-7,-7)(-8.48,-5.7)(-9.24,-3.83)
  \put(-8.7,3.3){\vector(-1,-2){0}}
}

\renewenvironment{proof}{
  \trivlist
  \item[\hskip\labelsep{\noindent\textbf{Proof.}}]
}{
  \qed
  \endtrivlist
}

\makeatother

\makeatletter
\newtheoremstyle{tw}{}{}{\itshape}{\parindent}{\scshape}{.}{ }{}
\newtheoremstyle{tw1}{}{}{\itshape}{\parindent}{\scshape}{.}{ }{}
\makeatother

\newtheorem{tw}{Theorem}[section]

\newtheorem{Notheorem}[tw]{\hspace{-.33em}}

\newtheorem{corollary}[tw]{Corollary}

\newtheorem{proposition}[tw]{Proposition}

\newtheorem{tw1}{Theorem}

\newenvironment{arablist}{
  \begin{list}{\textup{(\arabic{enumi})}}{
    \usecounter{enumi}
    \setlength{\topsep}{0pt}
    \setlength{\parsep}{0pt}
    \setlength{\itemsep}{0pt}
    \setlength{\leftmargin}{0pt}
    \settowidth{\labelwidth}{(m)}
    \setlength{\itemindent}{\parindent}
    \addtolength{\itemindent}{\labelwidth}
    \addtolength{\itemindent}{\labelsep}
  }
}{
  \end{list}
}

\newenvironment{alphlist}{
  \begin{list}{\textup{(\alph{enumi})}}{
    \usecounter{enumi}
    \setlength{\topsep}{0pt}
    \setlength{\parsep}{0pt}
    \setlength{\itemsep}{0pt}
    \setlength{\leftmargin}{0pt}
    \settowidth{\labelwidth}{(m)}
    \setlength{\itemindent}{\parindent}
    \addtolength{\itemindent}{\labelwidth}
    \addtolength{\itemindent}{\labelsep}
  }
}{
  \end{list}
}

\makeatletter

\def\ps@paper{

  \let\@mkboth\@gobbletwo

  \def\@evenhead{%

     \parbox{\textwidth}{

       \upshape\footnotesize{\bf\thepage}\hfill%

       {\bf B{ocian, Holm and Skowro\'{n}ski}}}%

     }

  \def\@oddhead{%

    \parbox{\textwidth}{%

    \upshape\footnotesize%

    {\bf Derived Equivalences of Nonstandard Domestic Selfinjective 
      Algebras}%

    \hfill{\bf\thepage}}%

  }

  \def\@oddfoot{}

  \let\@evenfoot\@oddfoot

}

\def\ps@first{

  \let\@mkboth\@gobbletwo

  \def\@evenhead{}%

  \def\@oddhead{}%

  \def\@oddfoot{\hfill\footnotesize[\thepage]\hfill}%

  \let\@evenfoot\@oddfoot%

}

\makeatother

\begin{document}

\thispagestyle{empty}

\author{\bf{Rafa{\l} Bocian$^{\mathrm{1}}$, Thorsten Holm$^{\mathrm{2}}$,
and Andrzej Skowro\'{n}ski$^{\mathrm{1,*}}$}\vspace{1.5cm}\\
$^{\mathrm{1}}$Faculty of Mathematics and Computer Science,\\
Nicolaus Copernicus University,\\
Toru\'{n}, Poland\vspace{0.5cm}\\
$^{\mathrm{2}}$Department of Pure Mathematics,\\
University of Leeds,\\
Leeds LS2 9JT, United Kingdom}
\title{{\bf Derived Equivalence Classification of~Nonstandard
Selfinjective Algebras of~Domestic Type\vspace{1cm}}}

\date{}

\maketitle

\renewcommand{\thefootnote}{}
\footnotetext{$^{\mathrm{*}}$Correspondence: Andrzej
Skowro\'{n}ski, Faculty of Mathematics and Computer Science,
Nicolaus Copernicus
University, Chopina 12/18, 87-100 Toru\'{n}, Poland;\\
E-mail: skowron@mat.uni.torun.pl.}

\begin{center}\parbox{14cm}{
\section*{\protect\centering ABSTRACT} \label{abstract}

We give a complete derived equivalence classification of all
nonstandard representation-infinite domestic selfinjective
algebras over an algebraically closed field. As a consequence, a
complete stable equivalence classification of these algebras is
obtained.\vspace{0.5cm}

\noindent \textit{Key Words}: Selfinjective algebra; Brauer graph
algebra; Derived equivalence; Stable equivalence; Domestic type.}
\end{center}

\newpage


\section{\protect\centering\hspace{-0.7cm}.\hspace{0.7cm}INTRODUCTION}
\label{introduction}
\mbox{ }

Throughout the paper $K$ will denote a fixed algebraically closed
field. By an algebra we mean a finite dimensional $K$-algebra
(associative, with an identity), which we shall assume (without
loss of generality) to be basic and connected. For an algebra $A$,
we denote by $\md A$ the category of finite dimensional left
$A$-modules, by $\ulmod A$ the stable category of $\md A$ (modulo
projectives), and by $\dualnosc^{b}\!\lo\md A\po$ the derived
category of bounded complexes of modules from $\md A$. Two
algebras $A$ and $B$ are said to be stably equivalent if the
stable module categories $\ulmod A$ and $\ulmod B$ are equivalent.
Moreover, two algebras $A$ and $B$ are said to be derived
equivalent if the derived categories $\dualnosc^{b}\!\lo\md A\po$
and $\dualnosc^{b}\!\lo\md B\po$ are equivalent as triangulated
categories. An algebra $A$ is called selfinjective if the
projective $A$-modules are injective. It is proved in Rickard
(1989b) that derived equivalent selfinjective algebras are also
stably equivalent.

From Drozd's Tame and Wild Theorem (Drozd, 1980) the class of
algebras may be divided into two disjoint classes. One class consists
of the tame algebras for which the indecomposable modules occur,
in each dimension $d$, in a finite number of discrete and a finite
number of one-parameter families. The second class is formed by
the wild algebras whose representation theory comprises the
representation theories of all finite dimensional $K$-algebras.
Accordingly, a classification of the indecomposable finite
dimensional modules is feasible only for the tame algebras.

One central problem of modern representation theory is the
determination of the derived equivalence classes of tame
selfinjective algebras. This has been done for the selfinjective
algebras of finite representation type (see Asashiba, 1999),
the algebras of
dihedral, semidihedral and quaternion type (including tame blocks
of group algebras) (see Holm, 1997, 1999), and the symmetric algebras of
tubular type (see Bia\l{}kowski, Holm and Skowro\'nski, 2003a,
2003b).

In this paper, we are concerned with the problem of the derived
equivalence classification of all representation-infinite tame
selfinjective algebras of domestic type. Recall that an algebra
has domestic (representation) type if there exists a common bound
(independent of the fixed dimension) for the numbers of
one-parameter families of indecomposable modules. The Morita
equivalence classification of these algebras splits into two
cases: the standard algebras, which admit simply connected Galois
coverings, and the remaining nonstandard ones. The class of
standard representation-infinite selfinjective algebras of
domestic type coincides with the class of selfinjective algebras
of Euclidean type, that is, the orbit algebras $\widehat{B}/G$ of
the repetitive algebras $\widehat{B}$ of tilted algebras $B$ of
Euclidean type with respect to actions of admissible infinite
cyclic automorphism groups $G$ (see Skowro\'nski, 1989). We refer
to Bocian and Skowro\'nski (2003, 2005a, 2005b) for the
classification of these algebras, and to Bocian, Holm and
Skowro\'nski (2004, 2005) for the derived equivalence
classification of the symmetric algebras of Euclidean type and the
one-parametric standard selfinjective algebras, respectively. It
has been proved recently (see Bocian and Skowro\'nski, 2005c;
Skowro\'nski, 2005, and Section \ref{roz2}) that the class of
nonstandard representation-infinite selfinjective algebras of
domestic type consists of (modified) Brauer graph algebras
$\Omega\lo T\po$ of Brauer graphs $T$ with one loop. The aim of
this paper is to give the derived equivalence classification and
the stable equivalence classification of these algebras. In order
to formulate our main result, consider the following family of
bound quiver algebras
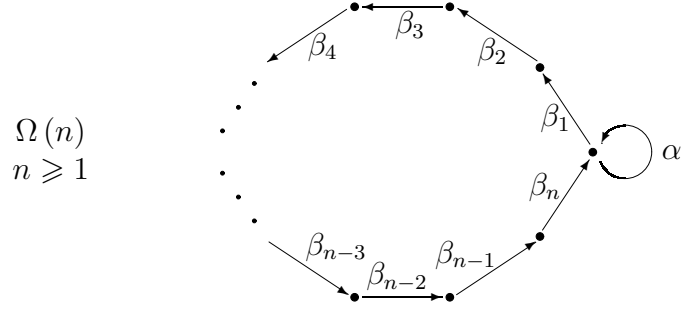
\begin{center}

\unitlength=1pt

\begin{picture}(311,176)(-145,-125)

 \put(99,-12){\vector(-2,3){17.8}}

 \put(79,19.5){\vector(-3,2){30}}

 \put(43,40){\vector(-1,0){30}}

 \put(7,39){\vector(-3,-2){30}}

 \put(81,-44.5){\vector(2,3){17.8}}

 \put(48.5,-69){\vector(3,2){30}}

 \put(13,-70){\vector(1,0){30}}

 \put(-22.5,-49){\vector(3,-2){30}}

 \put(80,17){\circle*{3}}

 \put(46,40){\circle*{3}}

 \put(10,40){\circle*{3}}

 \put(80,-47){\circle*{3}}

 \put(46,-70){\circle*{3}}

 \put(10,-70){\circle*{3}}

 \put(100,-15){\circle*{3}}

 \put(112,-15){\circlearrowr}

 \put(130,-15){\makebox(0,0){\normalsize $\alpha$}}

 \multiput(-28,11)(-6,-9){3}{\circle*{1.5}}

 \multiput(-40,-23)(6,-9){3}{\circle*{1.5}}

 \put(-105,-15){\makebox(0,0){\normalsize 
$\begin{array}{c} \Omega\lo n\po\\
 n \geqslant 1\end{array}$}}

 \put(86,-2){\makebox(0,0){\normalsize $\beta_{1}$}}
 \put(62,23){\makebox(0,0){\normalsize $\beta_{2}$}}

 \put(31,33){\makebox(0,0){\normalsize $\beta_{3}$}}

 \put(-2,25){\makebox(0,0){\normalsize $\beta_{4}$}}

 \put(82,-28){\makebox(0,0){\normalsize $\beta_{n}$}}

 \put(52,-54){\makebox(0,0){\normalsize $\beta_{n-1}$}}

 \put(26,-62){\makebox(0,0){\normalsize $\beta_{n-2}$}}

 \put(3,-51){\makebox(0,0){\normalsize $\beta_{n-3}$}}

 \put(10,-100){\makebox(0,0){\normalsize $\begin{array}{l}

 \alpha^{2}=\alpha\beta_{1}\beta_{2}\ldots\beta_{n},\;
\alpha\beta_{1}\beta_{2}\ldots\beta_{n}+\beta_{1}\beta_{2}
\ldots\beta_{n}\alpha=0,\\
 \beta_{n}\beta_{1}=0,\;\beta_{j}\beta_{j+1}\ldots\beta_{n}
\alpha\beta_{1}\beta_{2}\ldots\beta_{j-1}\beta_{j}=0,\;
 for\; 2\leqslant j \leqslant n.
 \end{array}$}}
\end{picture}
\end{center}

The following is the main result of this paper, providing a
complete derived equivalence classification and stable equivalence
classification of nonstandard representation-infinite
selfinjective algebras of domestic type.

\begin{tw1}\label{mainthm}

(1) Let $A=\Omega\lo T\po$, and assume that the Brauer graph $T$ has
$n$ edges. Then $A$ is derived equivalent (respectively, stably
equivalent) to $\Omega\lo n\po$.

(2) Any nonstandard representation-infinite selfinjective algebra
of domestic type is derived equivalent (resp. stably equivalent)
to an algebra $\Omega\lo n\po$.

Moreover, two algebras $\Omega\lo
m\po$ and $\Omega\lo n\po$ are derived equivalent (respectively,
stably equivalent) if and only if $m=n$.
\end{tw1}


\section{\protect\centering\hspace{-0.7cm}.\hspace{0.5cm} 
THE ALGEBRAS $\bf\Omega\lo T\po$}\label{roz2} \mbox{}

The aim of this section is to introduce the family $\Omega\lo
T\po$ of nonstandard Brauer graph algebras of domestic type.

A Brauer graph $T$ is a finite connected undirected graph, where
for each vertex there is a fixed circular order on the edges
adjacent to it. We draw $T$ in a plane and agree that the edges
adjacent to a given vertex are clockwise ordered. Here, we assume
that $T$ has exactly one cycle, which is a loop, being also its
direct successor. Therefore, $T$ is of the form
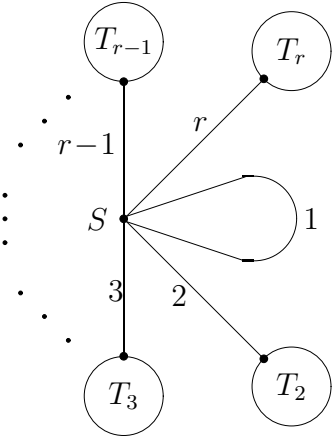
\begin{center}
\unitlength=1pt
\begin{picture}(151.00,165.00)(23.00,-23.00)
\put(81.00,59.00){\circle*{3.00}}
\put(81.00,111.00){\circle*{3.00}}
\put(81.00,7.00){\circle*{3.00}}
\put(134.00,112.00){\circle*{3.00}}
\put(134.00,6.00){\circle*{3.00}}
\put(81.00,60.00){\line(0,1){50.00}}
\put(82.00,60.00){\line(1,1){51.00}}
\put(81.00,58.00){\line(0,-1){50.00}}
\put(82.00,58.00){\line(1,-1){51.00}}
\put(152.00,59.00){\makebox(0,0){\normalsize $1$}}
\put(110.00,95.00){\makebox(0,0){\normalsize $r$}}
\put(67.00,87.00){\makebox(0,0){\normalsize $r\!-\!1$}}
\put(78.00,32.00){\makebox(0,0){\normalsize $3$}}
\put(102.00,30.00){\makebox(0,0){\normalsize $2$}}
\put(71.00,59.00){\makebox(0,0){\normalsize $S$}}
\multiput(60.00,105.00)(-9,-9){3}{\circle*{1.5}}
\multiput(60.00,13.00)(-9,9){3}{\circle*{1.5}}
\multiput(36.00,50.00)(0,9){3}{\circle*{1.5}}

\put(82.00,60.00){\line(3,1){45.00}}
\put(82.00,58.00){\line(3,-1){45.00}}
\put(126.00,59.00){\oval(40,32)[r]}
\put(81.00,126.00){\circle{30}} \put(81.00,-8.00){\circle{30}}
\put(144.50,122.50){\circle{30}} \put(144.50,-4.50){\circle{30}}
\put(81.00,126.00){\makebox(0,0){\normalsize $T_{r-1}$}}
\put(81.00,-8.00){\makebox(0,0){\normalsize $T_{3}$}}
\put(144.50,-4.50){\makebox(0,0){\normalsize $T_{2}$}}
\put(144.50,122.50){\makebox(0,0){\normalsize $T_{r}$}}
\end{picture}
\end{center}
where $T_{2}, T_{3},\ldots,T_{r-1},T_{r}$ are Brauer trees. A
Brauer graph $T$ defines a Brauer quiver $Q_{T}$ as follows:
\begin{alphlist}
\item the vertices of $Q_{T}$ correspond to the edges of $T$;
\item there is an arrow $i\longrightarrow j$ in $Q_{T}$ if and
only if $j$ is the consecutive edge of $i$ in the circular
ordering of the edges at a vertex of $T$.
\end{alphlist}
Observe that $Q_{T}$ is the union of (oriented) cycles, and every
vertex of $Q_{T}$ belongs to exactly two cycles. The cycles of
$Q_{T}$ are divided into two camps: $\alpha$-camps and
$\beta$-camps such that two cycles of $Q_{T}$ having nontrivial
intersection belong to different camps. We assume that the loop
corresponding to the unique loop of $T$ is denoted $\alpha_{1}$.
Moreover, the $\beta$-cycle of $Q_{T}$ corresponding to the vertex
$S$ of the loop of $T$ is said to be exceptional. For each vertex
$i$ of $Q_{T}$, we have
\begin{itemize}
\item $i\stackrel{\alpha_{i}}{\longrightarrow}\alpha\lo i\po$, the
arrow in the $\alpha$-camp of $Q_{T}$ starting at $i$;
\item $i\stackrel{\beta_{i}}{\longrightarrow}\beta\lo i\po$, the
arrow in the $\beta$-camp of $Q_{T}$ starting at $i$;
\end{itemize}
and the oriented cycles
\begin{displaymath}
A_{i}=\alpha_{i}\alpha_{\alpha\lo
i\po}\ldots\alpha_{\alpha^{-1}\!\lo i\po},\mbox{\hspace{0.5cm}}
B_{i}=\beta_{i}\beta_{\beta\lo i\po}\ldots\beta_{\beta^{-1}\!\lo
i\po}
\end{displaymath}
around the vertex $i$. Moreover, for each vertex $j$ of the
exceptional $\beta$-cycle different from $1$, consider the
oriented cycle
\begin{displaymath}
B_{i}^{\prime}=\beta_{j}\beta_{j+1}\ldots\beta_{r}\alpha_{1}\beta_{1}
\beta_{2}\ldots\beta_{j-1}.
\end{displaymath}
Finally, we define the algebra $\Omega\lo T\po$ as the bound
quiver algebra $KQ_{T}/J_{T}$, where $J_{T}$ is the ideal of the
path algebra $KQ_{T}$ generated by the elements:
\begin{arablist}
\item $\beta_{\beta^{-1}\lo i\po}\alpha_{i}$, $\alpha_{\alpha^{-1}\lo
i\po}\beta_{i}$, for all vertices $i$ of $Q_{T}$ different from
$1$,
\item $\beta_{r}\beta_{1}$,
\item $A_{i}-B_{i}$, for all vertices $i$ of $Q_{T}$ not lying on
the exceptional $\beta$-cycle,
\item $A_{j}-B^{\prime}_{j}$, for all vertices $j$ of the exceptional
$\beta$-cycle different from $1$,
\item $\alpha_{1}^{2}-\alpha_{1}\beta_{1}\ldots\beta_{r}$,
$\alpha_{1}\beta_{1}\ldots\beta_{r}+\beta_{1}\ldots\beta_{r}\alpha_{1}$.
\end{arablist}
\bigskip

We then have the following theorem (see Bocian and Skowro\'nski,
2005; Skowro\'nski, 2005).
\begin{tw}\label{tw21}
The Brauer graph algebras $\Omega\lo T\po$ form a complete family
of nonstandard representation-infinite selfinjective algebras of
domestic type.
\end{tw}

For each positive integer $n$, denote by $\Omega\lo n\po$ the
Brauer graph algebra $\Omega\lo T\lo n\po\po$ of the Brauer
loop-star $T\lo n\po$ of the form
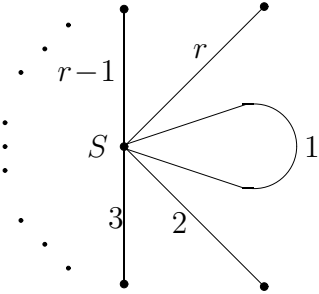
\begin{center}
\unitlength=1pt
\begin{picture}(121.00,115.00)(33.00,3.00)
\put(81.00,59.00){\circle*{3.00}}
\put(81.00,111.00){\circle*{3.00}}
\put(81.00,7.00){\circle*{3.00}}
\put(134.00,112.00){\circle*{3.00}}
\put(134.00,6.00){\circle*{3.00}}
\put(81.00,60.00){\line(0,1){50.00}}
\put(82.00,60.00){\line(1,1){51.00}}
\put(81.00,58.00){\line(0,-1){50.00}}
\put(82.00,58.00){\line(1,-1){51.00}}
\put(152.00,59.00){\makebox(0,0){\normalsize $1$}}
\put(110.00,95.00){\makebox(0,0){\normalsize $r$}}
\put(67.00,87.00){\makebox(0,0){\normalsize $r\!-\!1$}}
\put(78.00,32.00){\makebox(0,0){\normalsize $3$}}
\put(102.00,30.00){\makebox(0,0){\normalsize $2$}}
\put(71.00,59.00){\makebox(0,0){\normalsize $S$}}
\multiput(60.00,105.00)(-9,-9){3}{\circle*{1.5}}
\multiput(60.00,13.00)(-9,9){3}{\circle*{1.5}}
\multiput(36.00,50.00)(0,9){3}{\circle*{1.5}}

\put(82.00,60.00){\line(3,1){45.00}}
\put(82.00,58.00){\line(3,-1){45.00}}
\put(126.00,59.00){\oval(40,32)[r]}
\end{picture}
\end{center}

The following proposition describes the structure of the stable
Auslander-Reiten quivers $\Gamma_{\Omega\lo n\po}^{s}$ of the
algebras $\Omega\lo n\po$.
\begin{proposition}\label{stw22}
The stable Auslander-Reiten quiver $\Gamma_{\Omega\lo n\po}^{s}$
of $\Omega\lo n\po$ consists of a Euclidean component of type
$\mathbbm{Z}\widetilde{\mathbbm{A}}_{2n-1}$, a stable tube of rank
$2n\!-\! 1$, and a $K$-family of stable tubes of rank $1$.
\end{proposition}
\begin{proof}
Fix $n\geqslant 1$ and consider the bound quiver algebra $A\lo
n\po$
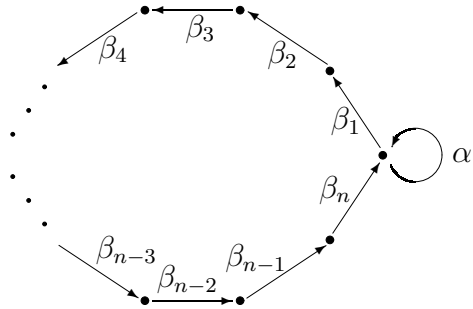
\begin{center}
\unitlength=1pt
\begin{picture}(281,176)(-90,-125)
 \put(99,-12){\vector(-2,3){17.8}}
 \put(79,19.5){\vector(-3,2){30}}
 \put(43,40){\vector(-1,0){30}}
 \put(7,39){\vector(-3,-2){30}}

 \put(81,-44.5){\vector(2,3){17.8}}
 \put(48.5,-69){\vector(3,2){30}}
 \put(13,-70){\vector(1,0){30}}
 \put(-22.5,-49){\vector(3,-2){30}}

 \put(80,17){\circle*{3}}
 \put(46,40){\circle*{3}}
 \put(10,40){\circle*{3}}

 \put(80,-47){\circle*{3}}
 \put(46,-70){\circle*{3}}
 \put(10,-70){\circle*{3}}

 \put(100,-15){\circle*{3}}

 \put(112,-15){\circlearrowr}
 \put(130,-15){\makebox(0,0){\normalsize $\alpha$}}

 \multiput(-28,11)(-6,-9){3}{\circle*{1.5}}
 \multiput(-40,-23)(6,-9){3}{\circle*{1.5}}

 \put(86,-2){\makebox(0,0){\normalsize $\beta_{1}$}}
 \put(62,23){\makebox(0,0){\normalsize $\beta_{2}$}}
 \put(31,33){\makebox(0,0){\normalsize $\beta_{3}$}}
 \put(-2,25){\makebox(0,0){\normalsize $\beta_{4}$}}

 \put(82,-28){\makebox(0,0){\normalsize $\beta_{n}$}}
 \put(52,-54){\makebox(0,0){\normalsize $\beta_{n-1}$}}
 \put(26,-62){\makebox(0,0){\normalsize $\beta_{n-2}$}}
 \put(3,-51){\makebox(0,0){\normalsize $\beta_{n-3}$}}

 \put(50,-100){\makebox(0,0){\normalsize $\begin{array}{l}
 \alpha^{2}=0,\;\beta_{n}\beta_{1}=0,\;\alpha\beta_{1}\beta_{2}\ldots\beta_{n}+\beta_{1}\beta_{2}\ldots\beta_{n}\alpha=0,\\
 \beta_{j}\beta_{j+1}\ldots\beta_{n}\alpha\beta_{1}\beta_{2}\ldots\beta_{j-1}\beta_{j}=0,\;
 for\; 2\leqslant j \leqslant n.
 \end{array}$}}
\end{picture}
\end{center}
Then the algebras $\Omega\lo n\po$ and $A\lo n\po$ are socle
equivalent, that is, the algebras $\Omega\lo n\po/\soc\Omega\lo
n\po$ and $A\lo n\po/\soc A\lo n\po$ are isomorphic, and
consequently the stable Auslander-Reiten quivers
$\Gamma_{\Omega\lo n\po}^{s}$ and $\Gamma_{A\lo n\po}^{s}$ are
isomorphic, by the known form
\begin{displaymath}
0\longrightarrow\rad P\longrightarrow\rad P/\soc P\oplus P
\longrightarrow P/\soc P\longrightarrow 0
\end{displaymath}
of Auslander-Reiten sequences with middle term having an
indecomposable projective-injective direct summand $P$ (see
Auslander, Reiten and Smalo, 1995). On the other hand, $A\lo n\po$
is a weakly symmetric standard selfinjective algebra
$\widehat{B\lo n\po}/\lo\sigma_{n}\varphi_{n}\po$ of Euclidean
type $\widetilde{\mathbbm{A}}_{2n-1}$, where $\widehat{B\lo n\po}$
is the repetitive algebra of the tilted algebra $B\lo n\po$ of
type $\widetilde{\mathbbm{A}}_{2n-1}$ given by the quiver
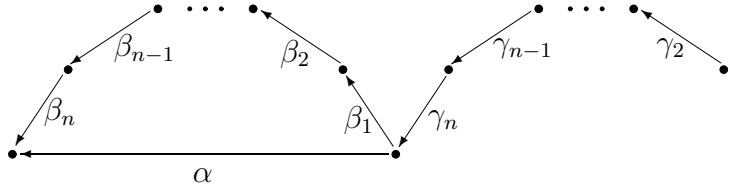
\begin{center}
\unitlength=1pt
\begin{picture}(278,80)(-48,-30)
 \put(99,-12){\vector(-2,3){17.8}}
 \put(119,14.5){\vector(-2,-3){17.8}}
 \put(79,19.5){\vector(-3,2){30}}
 \put(-26,15){\vector(-2,-3){17.8}}
 \put(7,39){\vector(-3,-2){30}}
 \put(97,-15){\vector(-1,0){139}}

 \put(80,17){\circle*{3}}
 \put(120,17){\circle*{3}}
 \put(46,40){\circle*{3}}
 \put(154,40){\circle*{3}}
 \put(10,40){\circle*{3}}
 \put(190,40){\circle*{3}}
 \put(-24,17){\circle*{3}}

 \put(-45,-15){\circle*{3}}
 \put(100,-15){\circle*{3}}

 \multiput(34,40)(-6,0){3}{\circle*{1.5}}
 \multiput(178,40)(-6,0){3}{\circle*{1.5}}

 \put(86,-2){\makebox(0,0){\normalsize $\beta_{1}$}}
 \put(118,-2){\makebox(0,0){\normalsize $\gamma_{n}$}}
 \put(148,25){\makebox(0,0){\normalsize $\gamma_{n-1}$}}
 \put(62,23){\makebox(0,0){\normalsize $\beta_{2}$}}
 \put(-27,1){\makebox(0,0){\normalsize $\beta_{n}$}}
 \put(204,25){\makebox(0,0){\normalsize $\gamma_{2}$}}
 \put(5,25){\makebox(0,0){\normalsize $\beta_{n-1}$}}
 \put(27,-23){\makebox(0,0){\normalsize $\alpha$}}

 \put(151,39){\vector(-3,-2){30}}
 \put(224,17){\circle*{3}}
 \put(223,19.5){\vector(-3,2){30}}
\end{picture}
\end{center}
and the relation $\gamma_{n}\beta_{1}=0$, $\varphi_{n}$ is the
canonical automorphism of $\widehat{B\lo n\po}$ whose square
$\varphi_{n}^{2}$ is the Nakayama automorphism $\nu_{\widehat{B\lo
n\po}}$ of $\widehat{B\lo n\po}$, and $\sigma_{n}$ is a rigid
automorphism of $\widehat{B\lo n\po}$ induced by the scalar
multiplication $\lo -1\po\alpha$ of $\alpha$. Then the required
shape of $\Gamma_{A\lo n\po}^{s}$, and hence of $\Gamma_{\Omega\lo
n\po}^{s}$, follows from Assem, Nehring and Skowro\'nski (1989)
and Skowro\'nski (1989).
\end{proof}

The following corollary is a direct consequence of the above
proposition.
\begin{corollary}\label{wn23}
Let $m$ and $n$ be different positive integers. Then the algebras
$\Omega\lo m\po$ and $\Omega\lo n\po$ are not stably equivalent.
\end{corollary}


\section{\protect\centering\hspace{-0.7cm}.\hspace{0.7cm}
PROOF OF THEOREM \ref{mainthm}}

Note that the part (2) of Theorem \ref{mainthm} follows from the
part (1) together with Theorem \ref{tw21}. Moreover, by (Rickard,
1989b), derived equivalent selfinjective algebras are stably
equivalent, and the final part of Theorem \ref{mainthm} follows
from Corollary \ref{wn23}.

So we have to prove that every algebra $\Omega(T)$ as defined in
Section \ref{roz2} is derived equivalent to the normal form
$\Omega(n)$, where $n$ is the number of edges of the Brauer graph
$T$.

Actually, we will give two different proofs of Theorem 1. The
first is based on Rickard's construction of tilting complexes for
Brauer tree algebras (see Rickard 1989b). We slightly have to
adapt the original construction for our purposes. This proof is
quite elegant, but for most technical details we will have to
refer to Rickard's paper. For the convenience of the reader who
might not be familiar with Rickard's paper we give a second, more
elementary and self-contained proof of Theorem 1. We give a
construction of an easy tilting complex whose endomorphism ring is
of the form $\Omega(T')$ where the exceptional cycle in $T'$ has
one vertex more than the cycle in $T$. Inductively, we can get all
vertices inside the cycle, up to derived equivalence, that is, we
get a normal form $\Omega(n)$.

The first proof will be given in subsection \ref{shrinking}, the second
proof in subsection \ref{enlarging}.

\subsection{Shrinking the Brauer trees} \label{shrinking}
Recall the definition of the algebras
$\Omega(T)$ from Section \ref{roz2}. The Brauer graph $T$ has the 
following shape
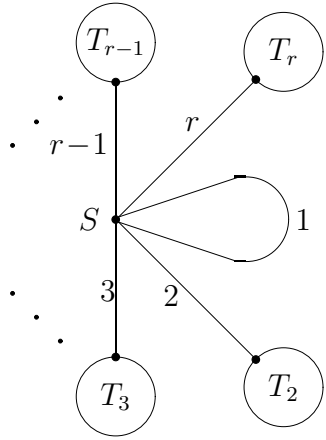
\begin{center}
\unitlength=1pt
\begin{picture}(151.00,165.00)(23.00,-23.00)
\put(81.00,59.00){\circle*{3.00}}
\put(81.00,111.00){\circle*{3.00}}
\put(81.00,7.00){\circle*{3.00}}
\put(134.00,112.00){\circle*{3.00}}
\put(134.00,6.00){\circle*{3.00}}
\put(81.00,60.00){\line(0,1){50.00}}
\put(82.00,60.00){\line(1,1){51.00}}
\put(81.00,58.00){\line(0,-1){50.00}}
\put(82.00,58.00){\line(1,-1){51.00}}
\put(152.00,59.00){\makebox(0,0){\normalsize $1$}}
\put(110.00,95.00){\makebox(0,0){\normalsize $r$}}
\put(67.00,87.00){\makebox(0,0){\normalsize $r\!-\!1$}}
\put(78.00,32.00){\makebox(0,0){\normalsize $3$}}
\put(102.00,30.00){\makebox(0,0){\normalsize $2$}}
\put(71.00,59.00){\makebox(0,0){\normalsize $S$}}
\multiput(60.00,105.00)(-9,-9){3}{\circle*{1.5}}
\multiput(60.00,13.00)(-9,9){3}{\circle*{1.5}}
\multiput(36.00,50.00)(0,9){3}{\circle*{1.5}}

\put(82.00,60.00){\line(3,1){45.00}}
\put(82.00,58.00){\line(3,-1){45.00}}
\put(126.00,59.00){\oval(40,32)[r]}
\put(81.00,126.00){\circle{30}} \put(81.00,-8.00){\circle{30}}
\put(144.50,122.50){\circle{30}} \put(144.50,-4.50){\circle{30}}
\put(81.00,126.00){\makebox(0,0){\normalsize $T_{r-1}$}}
\put(81.00,-8.00){\makebox(0,0){\normalsize $T_{3}$}}
\put(144.50,-4.50){\makebox(0,0){\normalsize $T_{2}$}}
\put(144.50,122.50){\makebox(0,0){\normalsize $T_{r}$}}
\end{picture}
\end{center}
where $T_{2}, T_{3},\ldots,T_{r-1},T_{r}$ are arbitrary Brauer
trees.

The construction below of a suitable tilting complex is inspired by
J. Rickard's construction of tilting complexes for Brauer tree algebras
(see Rickard, 1989b). Of course, in our situation we slightly have to
adapt this construction because we are dealing with Brauer graphs 
containing a cycle.

Recall that the edges of the Brauer graph $T$ correspond to the
vertices of the Brauer quiver $Q_T$, that is, to the simple
modules of the algebra $\Omega(T)$. For each edge $z$ of the
Brauer graph $T$ we shall define a bounded complex $Q(z)$ of
projective $\Omega(T)$-modules.

We then form the direct sum complex $Q:=\oplus_{z\in T} Q(z)$
over all edges of $T$
and conclude that it is actually a tilting complex for $\Omega(T)$.

For every edge $z$ corresponding to a vertex on the exceptional cycle of
$Q_T$ we set
$Q(z)$ to be the stalk complex with the projective indecomposable module
$P(z)$ concentrated in degree 0.

Now let $z$ be any edge not corresponding to a vertex on the
exceptional cycle. Then $z$ is contained in one of the Brauer
trees, say in $T_{i}$. Then there is a unique shortest path in $T$
from the edge $i$ to the edge $z$. Denote the edges on this path
by $z_0=i,z_1,\ldots,z_r=z$. From the usual Brauer tree relations,
we see that, up to scalar multiplication, there is a unique
homomorphism $P(z_j)\to P(z_{j+1})$ between the corresponding
projective indecomposable modules. Hence we get the following
complex of $\Omega(T)$-modules
$$Q(z)\,:\,0\to P(i)\to P(z_1)\to \ldots\to P(z_{r-1})\to P(z_r)\to 0$$
in which all maps are non-zero, and where $P(i)$ is in degree 0.

\begin{proposition} \label{tilt-T}
The complex $Q:= \oplus_{z\in T} Q(z)$ is a tilting complex for 
$\Omega(T)$.
\end{proposition}

\begin{proof}
Exactly the same as in (Rickard 1989b) for Brauer tree algebras.
In fact, the proof there carries over verbatim since there are no
non-zero homomorphisms $P(z)\to P(z')$ unless $z$ and $z'$ are
both corresponding to vertices on the exceptional cycle of $Q_T$,
or are both in the same Brauer tree $T_{i}$. We refrain from
reproducing the proof here and refer to (Rickard, 1989b) for
details.
\end{proof}

From Rickard's celebrated criterion (see Rickard, 1989a) we deduce
that the endomorphism ring (in the homotopy category) of the
tilting complex $Q$ is derived equivalent to $\Omega(T)$. The
following result describes this endomorphism ring $\End(Q)$ in the
homotopy category. In particular, the part (2) of the following
proposition completes the proof of Theorem 1.

Recall that the edges of $T$ on
the exceptional cycle
are denoted $1,2,\ldots,r$ (in the cyclic order). Whenever we write
$\Hom$ or $\End$ without index, we mean morphisms in the homotopy 
category of complexes.

\begin{proposition}
\begin{arablist}
\item The Cartan matrix of the endomorphism ring $\End(Q)$ has the
following entries $\tilde{c}_{z,z'}:=\dim\,\Hom(Q(z),Q(z'))$:
$$\tilde{c}_{z,z'} = \left\{
    \begin{array}{rcl}
       4 & ~~~ & \mbox{if $z=z'=1$} \\
       2 & ~~~ & \mbox{if $z=z'\neq 1$} \\
       2 & ~~~ & \mbox{if $1=z\neq z'$~or~$z\neq z'=1$} \\
       1 & ~~~ & \mbox{else}
    \end{array}
    \right.
$$
\item $\End(Q)$ is isomorphic to the normal form $\Omega(n)$, 
where $n$ is the number of edges of $T$.
\end{arablist}
\end{proposition}

\begin{proof}
(1) The Cartan invariants of the endomorphism ring of a tilting
complex can conveniently be computed from the Cartan matrix of
$\Omega(T)$ using a well-known alternating sum formula (see
Happel, 1988): let $Q(z)=(Q_r)_{r\in\mathbb{Z}}$ and
$Q(z')=(Q'_s)_{s\in\mathbb{Z}}$, then
$$\tilde{c}_{z,z'} = \sum_{r,s} (-1)^{r-s}
\dim\,\Hom_{\Omega(T)}(Q_r,Q_s').
$$
(Note that the sum is indeed finite since we are dealing with
bounded complexes.)

In our situation these computations are particularly easy. In
fact, from the relations of $\Omega(T)$ we see that
$\Hom_{\Omega(T)}(P(v),P(v'))=0$ unless the edges $v,v'$ of $T$
(resp. the vertices of $Q_T$) both belong to the same simple
cycle.

We then leave the details of the slightly tedious, but straightforward
computations of the Cartan invariants $\tilde{c}_{z,z'}$ to the reader.
\smallskip

(2) Note that by part (1) the endomorphism ring has exactly the same
Cartan matrix as $\Omega(n)$. So it suffices to define homomorphisms
between the summands of $Q$, corresponding to the arrows of
$\Omega(n)$, and to show that these maps satisfy the defining relations
of $\Omega(n)$, up to homotopy.

Consider first the summands $Q(z)$ for all $z$ in a fixed Brauer
tree $T_{i}$ (including the vertex $i$ on the cycle). Then the
above tilting complex construction is precisely Rickard's
construction and we know the endomorphism ring $\End(\oplus_{z\in
B_i} Q(z))$ from (Rickard, 1989b). It is again a Brauer tree,
where the tree is a star (without multiplicity); the maps
corresponding to the arrows of the Brauer quiver are in degree 0
given as follows (we don't need to know the maps in the other
degrees precisely). Let $Q(z_i)$ be the direct successor of $Q(i)$
in the ordering of the star. The map $Q(i)\to Q(z_i)$ is given by
multiplication with $B_i$ (the $\beta$-cycle at $i$); all arrows
starting in $Q(z)$, where $z\neq i$, are given by the identity on
$P(i)$.

Now we consider the 'global' picture involving all trees $T_{i}$.
The crucial observation is that the map $Q(i)\to Q(z_i)$ factors
through all the other summands of $Q$. In fact, for any $i$,
define homomorphisms $Q(i)\to Q(z_{i+1})$ by multiplication with
$\beta_i$ in degree 0. (Note that this is a homomorphism of
complexes since $\beta\alpha=0$ in $\Omega(T)$.)

So we have a cyclic ordering of all summands of $Q$ as follows
(instead of $Q(z)$ we just write $z$ for abbreviation, and
'vertices of $T_{i}$' means vertices $\neq i,z_i$, the ordering of
them as in the star they form by Rickard's argument):
$$1,z_2,(\mbox{vertices of $T_{2}$}), 2, z_3, (\mbox{vertices of
$T_{3}$}),\ldots, r-1, z_{r}, (\mbox{vertices of $T_{r}$}), r,
1.$$ The maps $Q(i)\to Q(z_{i+1})$ are as defined above, and maps
between the $Q(z)$, $z\in B_i$ as in Rickard's construction, that
is, in degree 0 given by the identity.

With these definitions, it is straightforward to check that
$\End(Q)$ satisfies the defining relations of $\Omega(n)$. In
fact, in the cyclic ordering described above, most maps are in
degree 0 the identity, the remaining ones are given by
multiplication with the arrows $\beta_i$ on the exceptional cycle
of $T$. So the desired relations of $\End(Q)$ follow directly from
the relations of $\Omega(T)$ involving the exceptional cycle.

This completes the proof of the proposition.
\end{proof}


\subsection{Enlarging the cycle} \label{enlarging}

We now give our second, more self-contained, inductive
proof of Theorem 1. In each step we will show that we
can enlarge the exceptional cycle, up to derived equivalence.
Inductively, we can get all vertices inside the cycle and
obtain one of the normal forms $\Omega(n)$.

Recall the definition of $\Omega(T)$ and the shape of the Brauer 
graph $T$
from Section \ref{roz2}. Also recall that in the corresponding Brauer
quiver $Q_T$ the exceptional cycle was set to be a $\beta$-cycle.

Now consider the vertex $2$ in $Q_T$ (corresponding to the edge
$2$ in $T$). We can assume that the Brauer tree $T_{2}$ is not
empty (otherwise take the first vertex on the cycle with non-empty
Brauer tree).

Let $2_1$ denote the direct successor of $2$ on the $\alpha$-cycle
at $2$. Moreover, let $2^1,\ldots,2^k$
be the vertices on the $\beta$-cycle attached to $2_1$.
(If existing, otherwise $P(2^k)$ does not occur in the
tilting complex to be defined below.)

We define the following complexes of projective $\Omega(T)$-modules.
Set
$$Q'(2_1)\,:\,0\to P(2)\oplus P(2^k)
\stackrel{(\alpha,\beta)}{\longrightarrow}
P(2_1)\to 0$$
where $P(2)$ is in degree 0.
For all other vertices $z\neq 2_1$ in $Q_{T}$ we let
$Q'(z)$ be the stalk complex with $P(z)$ concentrated in degree 0.

\begin{proposition} \label{tilt-T2}
The complex $Q':=\oplus_{z\in T} Q'(z)$
is a tilting complex for $\Omega(T)$.
\end{proposition}

\begin{proof}
From the definition of tilting complex, there are two properties
to check:

(i) $\add(Q)$ generates the homotopy category of bounded complexes
of projective $\Omega(T)$-modules as triangulated category,

(ii) $\Hom(Q,Q[r])=0$ for all $r\neq 0$ (up to homotopy).
\smallskip

For (i), just observe that $P(2_1)$ is homotopy equivalent
to the mapping cone of
the obvious map $Q(2_1)\to Q(2)\oplus Q(2^k)$ given by the identity
in degree 0.

For (ii), first observe that $\Hom(Q,Q[r])=0$ for all $|r|\ge 2$,
since we are dealing with two-term complexes.

Let us consider the case $r=1$. Every homomorphism $P(z)\to
P(2_1)$, where $z\neq 2_1$, factors through $P(2)\oplus
P(2^k)\stackrel{(\alpha,\beta)}{\longrightarrow} P(2_1)$ (in fact,
every path in $Q_{T}$ going to the vertex $2_1$ ends with the
arrow $\alpha$ from $2$ to $2_1$ or with the arrow $\beta$ from
$2^k$ to $2_1$). Hence, $\Hom(Q',Q'[1])=0$ up to homotopy.

Now consider $r=-1$. It follows from the (Brauer tree) relations
of $\Omega(T)$ that no non-zero homomorphism $P(2_1)\to P(z)$,
where $z\neq 2_1$, gives the zero map when composed with right
multiplication by $\alpha:P(2)\to P(2_1)$ and with
$\beta:P(2^k)\to P(2_1)$. Hence, $\Hom(Q',Q'[-1])=0$.
\end{proof}

As a consequence, $\Omega(T)$ is derived equivalent to the
endomorphism ring of the tilting complex $Q'$.
The latter is described in detail by the following result.
The crucial aspect is that $\End(Q')$ is of the form
$\Omega(T')$ for some Brauer graph $T'$ whose exceptional
cycle has one vertex more than the one for $T$.

\begin{proposition} \label{end-T2}
The endomorphism ring (in the homotopy category)
of the tilting complex $Q'$ is an algebra $\Omega(T')$ where the
Brauer graph $T'$ is of the following shape
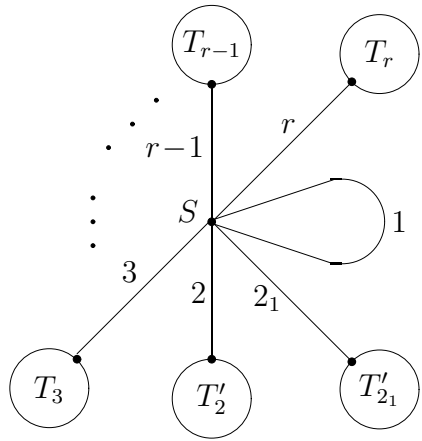
\begin{center}
\unitlength=1pt
\begin{picture}(151.00,165.00)(23.00,-23.00)
\put(81.00,59.00){\circle*{3.00}}
\put(81.00,111.00){\circle*{3.00}}
\put(81.00,7.00){\circle*{3.00}}
\put(134.00,112.00){\circle*{3.00}}
\put(134.00,6.00){\circle*{3.00}}
\put(30.00,7.00){\circle*{3.00}}

\put(81.00,60.00){\line(0,1){50.00}}
\put(81.00,60.00){\line(-1,-1){51.00}}
\put(82.00,60.00){\line(1,1){51.00}}
\put(81.00,58.00){\line(0,-1){50.00}}
\put(82.00,58.00){\line(1,-1){51.00}}
\put(152.00,59.00){\makebox(0,0){\normalsize $1$}}
\put(110.00,95.00){\makebox(0,0){\normalsize $r$}}
\put(67.00,87.00){\makebox(0,0){\normalsize $r\!-\!1$}}
\put(76.00,32.00){\makebox(0,0){\normalsize $2$}}
\put(50.00,40.00){\makebox(0,0){\normalsize $3$}}
\put(102.00,30.00){\makebox(0,0){\normalsize $2_1$}}
\put(72.00,63.00){\makebox(0,0){\normalsize $S$}}
\multiput(60.00,105.00)(-9,-9){3}{\circle*{1.5}}
\multiput(36.00,50.00)(0,9){3}{\circle*{1.5}}

\put(82.00,60.00){\line(3,1){45.00}}
\put(82.00,58.00){\line(3,-1){45.00}}
\put(126.00,59.00){\oval(40,32)[r]}
\put(81.00,126.00){\circle{30}}
\put(81.00,-8.00){\circle{30}}
\put(144.50,122.50){\circle{30}}
\put(144.50,-4.50){\circle{30}}
\put(19.50,-4.00){\circle{30}}

\put(81.00,126.00){\makebox(0,0){\normalsize $T_{r-1}$}}
\put(81.00,-8.00){\makebox(0,0){\normalsize $T_{2}'$}}
\put(144.50,-4.50){\makebox(0,0){\normalsize $T_{2_1}'$}}
\put(19.50,-4.50){\makebox(0,0){\normalsize $T_{3}$}}
\put(144.50,122.50){\makebox(0,0){\normalsize $T_{r}$}}
\end{picture}
\end{center}
where $T_2'$ and $T_{2_1}'$ are Brauer trees with less edges than
$T_2$. All other Brauer trees $T_3,\ldots,T_r$ remain unchanged.

In particular, $\Omega(T)$ is derived equivalent to
$\Omega(T')$.
\end{proposition}

\begin{proof} Note that apart from $Q(2_1)$, all summands in the
tilting complex $Q'$ are stalk complexes. So it is clear that
$\End(Q')$ when compared with $\Omega(T)$ only changes around
the vertex $2_1$. We will indicate below the major changes,
leaving some details of the straightforward
verifications to the reader. We freely use the notation introduced
in the proof of the previous proposition.

The crucial observation is that the map $\beta_1:Q(1)\to Q(2)$
between stalk complexes now factors through $Q(2_1)$. In fact,
it is the composition of the maps given in degree 0 by
$(\beta_1,0):Q(1)\to Q(2_1)$ and $(\id,0):Q(2_1)\to Q(2)$.
(Note that the first map is indeed a homomorphism of complexes
since $\beta_1\alpha=0$ in $\Omega(T)$.)
This explains that $2_1$ now belongs to the exceptional cycle
of the Brauer graph $T'$ of $\End(Q')$. (Clearly, the cycle
relations are not affected since the above composition is
just the 'old' map $\beta_1$.)

It remains to explain the structure of the new Brauer trees $T_2'$
and $T_{2_1}'$. Again, it suffices to describe the $\alpha$ and
$\beta$-cycles attached to the vertices $2$, $2_1$ and $2^k$; the
others remain unchanged, since the corresponding complexes are
stalk complexes.

Let $2_1,2_2,\ldots,2_j$ be the vertices on the 'old'
$\alpha$-cycle attached to $2$. Then the 'new' $\alpha$-cycle
attached to $2$ in $T_2'$ has vertices $2_2,\ldots,2_j$, and the
first map $Q(2)\to Q(2_2)$ is given by multiplication with
$\alpha^2$.

Otherwise, the Brauer tree $T_2'$ is the same as $T_2$. (But
$T_2'$ has now less vertices, for example, $2_1$ is missing.)

The $\alpha$-cycle at $2_1$ in $T'$ consists of the
following vertices: $2_1,2^k$ and the vertices of the
$\alpha$-cycle at $2^k$ in $T$. The first map
$Q(2_1)\to Q(2^k)$ is given by the projection in degree 0,
the last map back to $Q(2_1)$ is given by $(0,\alpha)$
in degree 0. (Note that this is indeed a homomorphism
of complexes by the usual Brauer tree relations.)

The 'new' $\beta$-cycle at $2^k$ has vertices $2^k,2^1,\ldots,
2^{k-1}$. Note that the 'old' vertex $2_1$ is no longer part
of this cycle. The first map $Q(2^k)\to Q(2^1)$ is now given
by $\beta^2$.

From the above descriptions it is clear that all the defining
relations of $\Omega(T')$ follow from the corresponding
relations of $\Omega(T)$. In fact, even if some cycles
in the quiver $Q_{T'}$ changed and became shorter,
the maps obtained by walking
around the cycle are exactly the same as before.

This completes the proof of the proposition.
\end{proof}

Now we are in the position to complete our second proof of
Theorem \ref{mainthm}.
\bigskip

\noindent
{\bf Proof of Theorem \ref{mainthm}.} Let $\Omega(T)$ be an algebra for
an arbitrary Brauer graph as defined in Section \ref{roz2}.
Now Proposition \ref{end-T2}
states that $\Omega(T)$ is derived equivalent
to an algebra $\Omega(T')$ where the Brauer graph $T'$ has one
edge more inside the cycle around the vertex $S$.

Inductively applying Proposition \ref{end-T2} (with vertex 2
replaced by any vertex on the cycle), we can move all edges inside
the cycle around $S$, up to derived equivalence. The resulting
algebra corresponds to a Brauer graph in which all attached Brauer
trees $T_j$ are empty. Clearly, this is the normal form algebra
$\Omega(n)$ where $n$ is the number of edges of the original
Brauer graph $T$.

This completes the proof of Theorem \ref{mainthm}. $\hfill\Box$


\section*{\protect \centering ACKNOWLEDGMENTS}
\mbox{ }

The first and third named authors acknowledge support from the
Polish Scientific Grant KBN No. 1 P03A 018 27.

\section*{\protect \centering REFERENCES}\mbox{}

Asashiba, H. (1999). The derived equivalence classification of
representation-finite selfinjective algebras. \textit{J. Algebra}
214:182--221.

Assem, I., Nehring, J., Skowro\'{n}ski, A. (1989). Domestic
trivial extensions of simply connected algebras. \textit{Tsukuba
J. Math.} 13:31--72.

Auslander, M., Reiten, I., Smalo, S. O. (1995). Representation
Theory of Artin Algebras. \textit{Cambridge Studies in Advanced
Math.} Vol. 36. Cambridge University Press.

Bia\l{}kowski, J., Holm, T., Skowro\'{n}ski, A. (2003a). Derived
equivalences for tame weakly symmetric algebras having only
periodic modules. \textit{J. Algebra} 269:652--668.

Bia\l{}kowski, J., Holm, T., Skowro\'{n}ski, A. (2003b). On
nonstandard tame selfinjective algebras having only periodic
modules. \textit{Colloq. Math.} 97:33--47.

Bocian, R., Holm, T., Skowro\'{n}ski, A. (2004). Derived
equivalence classification of weakly symmetric algebras of
Euclidean type. \textit{J. Pure Appl. Algebra} 191:43--74.

Bocian, R., Holm, T., Skowro\'{n}ski, A. (2005). Derived
equivalence classification of one-parametric selfinjective
algebras. \textit{Preprint, Toru\'n}.

Bocian, R., Skowro\'{n}ski, A. (2003). Symmetric special biserial
algebras of Euclidean type. \textit{Colloq. Math.} 96:121--148.

Bocian, R., Skowro\'{n}ski, A. (2005a). Weakly symmetric algebras
of Euclidean type. \textit{J. reine angew. Math.} 580:157--199.

Bocian, R., Skowro\'{n}ski, A. (2005b). One-parametric
selfinjective algebras. \textit{J .Math. Soc. Japan} 57:491--512.

Bocian, R., Skowro\'{n}ski, A. (2005c). Socle deformations of
selfinjective algebras of Euclidean type. \textit{Preprint,
Toru\'n}.

Drozd, Yu. A. (1980). Tame and wild matrix problems. \textit{In:
Representation Theory II.} Lecture Notes in Math. Vol. 832.
Berlin-New York: Springer-Verlag, pp. 242--258.

Happel, D. (1988). Triangulated Categories in Representation
Theory of Finite Dimensional Algebras. \textit{London Mathematical
Society.} Lecture Note Series. Vol. 119. Cambridge University
Press.

Holm, T. (1997). Derived equivalent tame blocks. \textit{J. Algebra}
194:178--200

Holm, T. (1999). Derived equivalence classification of algebras of
dihedral, semidihedral, and quaternion type. \textit{J. Algebra}
211:159--205.

Rickard, J. (1989a). Morita theory for derived categories.
\textit{J. London Math. Soc.} 39:436--456.

Rickard, J. (1989b). Derived categories and stable equivalence.
\textit{J. Pure Appl. Algebra} 61:303--317.

Skowro\'{n}ski, A. (1989). Selfinjective algebras of polynomial
growth. \textit{Math. Annalen.} 285:177--199.

Skowro\'{n}ski, A. (2005). Nonstandard selfinjective algebras of
domestic type. \textit{Preprint, Toru\'n}.

\end{document}